\documentclass[12pt,leqno]{article}

\usepackage{amsfonts}
\usepackage{amsmath}
\usepackage{amssymb}
\usepackage{epsf}

\newcommand{\Xcomment}[1]{}

\makeatletter
\renewcommand{\section}{\@startsection{section}{1}{0pt}%
{-3.5ex plus -1ex minus -.2ex}{2.3ex plus .2ex}%
{\large\bf}}
\makeatother

\newtheorem{theorem}{Theorem}[section]
\newtheorem{lemma}[theorem]{Lemma}
\newtheorem{corollary}[theorem]{Corollary}
\newtheorem{prop}[theorem]{Proposition}

\newenvironment{proof}{\noindent{\bf Proof.\/}}%
{$\qed$\vskip 0.1in}

\newenvironment{numitem}{\refstepcounter{equation}\begin{enumerate}%
\item[(\arabic{equation})]$\quad$}{\end{enumerate}}

\newcommand{\refeq}[1]{(\ref{eq:#1})}  

\def\qed{ \ \vrule width.2cm height.2cm depth0cm}

\def\tilde{\widetilde}

\def\bar{\overline}
\def\eps{\epsilon}

\def\Rset{{\mathbb R}}
\def\Zset{{\mathbb Z}}

\def\Qset{{\mathbb Q}}

\def\Cscr{\mathcal{C}}

\def\Lscr{\mathcal{L}}
\def\Mscr{\mathcal{M}}
\def\Pscr{\mathcal{P}}

\def\Rscr{\mathcal{R}}

\def\Wscr{\mathcal{W}}

\def\Pin{P^{\rm in}}
\def\Pout{P^{\rm out}}
\def\ellin{t}
\def\ellout{h}

\def\sprod{\;\hbox{\unitlength=1mm\begin{picture}(4,4)%
\put(0,0){$\triangleright$}\put(1.5,0){$\triangleleft$}%
\end{picture}}\;}

\begin{document}

\begin{center}
{\Large\bf Combinatorics of $A_2$-crystals}%
\footnote[1]{This research was supported in part by NWO--RFBR grant
047.011.2004.017 and by RFBR grant 05-01-02805
CNRSL\_a.}
 \end{center}

\begin{center}
{\sc Vladimir~I.~Danilov}\footnote[2]
{Central Institute of Economics and Mathematics
of the RAS, 47, Nakhimovskii Prospect, 117418 Moscow, Russia;
emails: danilov@cemi.rssi.ru, koshevoy@cemi.rssi.ru.},
{\sc Alexander~V.~Karzanov}\footnote[3]{Institute for System
Analysis of the RAS, 9, Prospect 60 Let Oktyabrya, 117312 Moscow,
Russia; email: sasha@cs.isa.ru. },\\
{\sc and Gleb~A.~Koshevoy}$^2$
\end{center}

\begin{quote}
{\bf Abstract.}
We show that a connected regular $A_2$-crystal (the crystal
graph of an irreducible representation of $sl_3$) can be produced
from two half-grids by replicating them and glying together in a
certain way. Also some extensions and related aspects are discussed.
  \end{quote}



\section{Introduction}  \label{sec:intr}

The notion of crystals introduced by Kashiwara~\cite{kas-90,kas-95}
embraces a wide class of edge-colored digraphs (directed graphs).
Stembridge~\cite{ste-03} pointed out a list of graph-theoretic axioms
characterizing the crystals with finite monochromatical paths that
are related to representations of simply-laced Lie algebras (i.e.,
with a Cartan matrix whose off-diagonal entries are 0 or --1), called
regular {\em simply-laced crystals}. Each of these axioms imposes a
simple local condition on a 2-colored subgraph of the digraph. In
particular, a digraph is a regular simple-laced crystal if and only
if each (inclusionwise) maximal 2-colored subgraph in it is
a regular simple-laced crystal (for a general
result of this type on any crystals of representations that have a
unique maximal vertex, see~\cite{KKM-92}). This shows an importance
of a proper study of 2-colored crystals.

This paper is the first in our series of works devoted to a
combinatorial study of crystals of representations, and to related
topics. Here we consider regular 2-colored simply-laced crystals
for the Cartan matrix $A_2=\binom{\;\;2\;-1}{-1\;\;2}$. When such a
crystal is connected and all monochromatic paths in it are finite,
we refer to it as an {\em RC-graph} (abbreviating ``regular crystal
graph''). Our main structural theorem says that an RC-graph
can be produced, by use replicating and glying together in a certain
way, from two RC-graphs of a very special form, viewed as triangular
halves of two-dimensional square grids.
As a result, the combinatorial structure of these objects becomes
rather transparent, giving rise to revealing additional properties of
RC-graphs and their extensions. In particular, it follows that an
RC-graph is the Hasse diagram of a finite lattice.

The paper is organized as follows. Section~\ref{sec:cryst} gives
basic definitions and exhibits some elementary properties of
RC-graphs. Section~\ref{sec:theo} contains
the formulation and proof of the above-mentioned structural
theorem. Section~\ref{sec:infin} explains how to
extend this result to 2-colored digraphs having infinite
monochromatic paths. Section~\ref{sec:patt}
describes natural embeddings of RC-graphs in the
Abelian groups $\Zset^4$ and $\Zset^3$. In particular, it explicitly
describes a relation to Gelfand-Tsetlin patterns.
The concluding Section~\ref{sec:univ} considers the disjoint union of
all RC-graphs, called the {\em universal RC-graph}, and
associate to it a certain semigroup of the Abelian group $\Zset^5$.

\section{RC-graphs} \label{sec:cryst}

Let $K$ be a digraph with vertex set $V$
and with edge set $E$ partitioned into two subsets $E_1,E_2$.
We say that an edge in $E_i$ has {\em color} $i$ and call it an
edge with {\em color} $i$, or, briefly, an $i$-{\em edge}.
Unless explicitly stated otherwise, any digraph in
question is assumed to be (weakly) connected, i.e., it is not
representable as the disjoint union of two nonempty digraphs.
An {\em RC-graph} is defined by
imposing on $K=(V,E_1,E_2)$ four axioms (A1)--(A4) described below.
(In fact, we reformulate axioms (P1)--(P6),(P5'),(P6')
given in~\cite{ste-03} for $n$-colored simply-laced crystals, to
our case $n=2$.)

The first axiom concerns the
structure of monochromatic subgraphs $(V,E_i)$.
\begin{itemize}
\item[(A1)] For $i=1,2$, each connected component of
$(V,E_i)$ is a finite simple (directed) path, i.e., a sequence
$(v_0,e_1,v_1,\ldots,e_k,v_k)$, where $v_0,v_1,\ldots,v_k$ are
distinct vertices and each $e_i$ is an edge going from $v_{i-1}$ to
$v_i$.
  \end{itemize}

In particular, each vertex $v$ has at most one outgoing 1-edge and at
most one incoming 1-edge, and similarly for 2-edges.
For convenience, we refer to a {\em maximal} monochromatic path in
$K$, with color $i$ on the edges, as an $i$-{\em line}.
The $i$-line passing through a given
vertex $v$ (possibly consisting of the only vertex $v$) is denoted by
$P_i(v)$, its part from the first vertex to $v$ by $\Pin_i(v)$, and
its part from $v$ to the last vertex by $\Pout_i(v)$.
The lengths of $\Pin_i(v)$ and of $\Pout_i(v)$ (i.e., the numbers of
edges in these paths) are denoted by $\ellin_i(v)$ and
$\ellout_i(v)$, respectively.

The second axiom tells us how these
lengths can change when one traverses an edge of the other color.
\begin{itemize}
\item[(A2)] Each $i$-line $P$ ($i=1,2$) contains a vertex $r$
satisfying the following property: for any edge $(u,v)$
(from a vertex $u$ to a vertex $v$) in $\Pin_i(r)$, one holds
$\ellin_{3-i}(v)=\ellin_{3-i}(u)-1$ and
$\ellout_{3-i}(v)=\ellout_{3-i}(u)$, and for any edge $(u',v')$ in
$\Pout_i(r)$, one holds $\ellin_{3-i}(v')=\ellin_{3-i}(u')$ and
$\ellout_{3-i}(v')=\ellout_{3-i}(u')+1$.
  \end{itemize}

\begin{figure}[htb]                  
 \begin{center}
  \unitlength=1mm
  \begin{picture}(90,50)
\put(0,5){\circle{1.0}}
\put(10,5){\circle{1.0}}
\put(20,5){\circle{1.0}}
\put(30,5){\circle{1.0}}
\put(40,5){\circle{1.0}}
\put(50,5){\circle{1.0}}
\put(10,13){\circle{1.0}}
\put(20,13){\circle{1.0}}
\put(30,13){\circle{1.0}}
\put(40,13){\circle{1.0}}
\put(50,13){\circle{1.0}}
\put(20,21){\circle{1.0}}
\put(30,21){\circle{1.0}}
\put(40,21){\circle{1.0}}
\put(50,21){\circle{1.0}}
\put(20,29){\circle{1.0}}
\put(30,29){\circle{1.0}}
\put(40,29){\circle{1.0}}
\put(50,29){\circle{1.0}}
\put(60,29){\circle{1.0}}
\put(20,37){\circle{1.0}}
\put(30,37){\circle{1.0}}
\put(40,37){\circle{1.0}}
\put(50,37){\circle{1.0}}
\put(60,37){\circle{1.0}}
\put(70,37){\circle{1.0}}
\put(20,45){\circle{1.0}}
\put(30,45){\circle{1.0}}
\put(40,45){\circle{1.0}}
\put(50,45){\circle{1.0}}
\put(60,45){\circle{1.0}}
\put(70,45){\circle{1.0}}
\put(80,45){\circle{1.0}}
\put(40,0){\vector(0,1){4.5}}
\put(40,5){\vector(0,1){7.5}}
\put(40,13){\vector(0,1){7.5}}
\put(40,21){\vector(0,1){7.5}}
\put(40,29){\vector(0,1){7.5}}
\put(40,37){\vector(0,1){7.5}}
\put(40,45){\line(0,1){5}}
\put(0,5){\vector(1,0){9.5}}
\put(10,5){\vector(1,0){9.5}}
\put(20,5){\vector(1,0){9.5}}
\put(30,5){\vector(1,0){9.5}}
\put(40,5){\vector(1,0){9.5}}
\put(10,13){\vector(1,0){9.5}}
\put(20,13){\vector(1,0){9.5}}
\put(30,13){\vector(1,0){9.5}}
\put(40,13){\vector(1,0){9.5}}
\put(20,21){\vector(1,0){9.5}}
\put(30,21){\vector(1,0){9.5}}
\put(40,21){\vector(1,0){9.5}}
\put(20,29){\vector(1,0){9.5}}
\put(30,29){\vector(1,0){9.5}}
\put(40,29){\vector(1,0){9.5}}
\put(50,29){\vector(1,0){9.5}}
\put(20,37){\vector(1,0){9.5}}
\put(30,37){\vector(1,0){9.5}}
\put(40,37){\vector(1,0){9.5}}
\put(50,37){\vector(1,0){9.5}}
\put(60,37){\vector(1,0){9.5}}
\put(20,45){\vector(1,0){9.5}}
\put(30,45){\vector(1,0){9.5}}
\put(40,45){\vector(1,0){9.5}}
\put(50,45){\vector(1,0){9.5}}
\put(60,45){\vector(1,0){9.5}}
\put(70,45){\vector(1,0){9.5}}
\put(41.5,47){P}
\put(40,21){\circle*{2}}
  \end{picture}
 \end{center}
 \caption{A example of changing the lengths $\ellin_1$ and $\ellout_1$
along line $P$ with color 2. The thick dot indicates the critical
vertex in $P$.}
  \label{fig:A2}
  \end{figure}
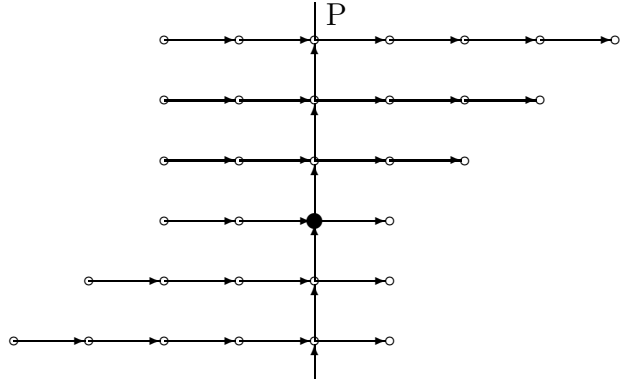

Clearly such a vertex $r$ is unique; it is called the {\em critical}
vertex of the given line $P$. Axiom (A2) is illustrated in
Fig.~\ref{fig:A2}.

The digraphs defined by axioms (A1),(A2) constitute a subclass of the
class of so-called (locally finite) {\em normal $A_2$-crystals} of
Kashiwara. (In the definition for the latter, axiom (A2) is replaced
by a relaxed axiom: for an $i$-edge $(u,v)$, either
$\ellin_{3-i}(v)=\ellin_{3-i}(u)-1$ and
$\ellout_{3-i}(v)=\ellout_{3-i}(u)$, or
$\ellin_{3-i}(v)=\ellin_{3-i}(u)$ and
$\ellout_{3-i}(v)=\ellout_{3-i}(u)+1$.)
We will refer to a digraph satisfying (A1) and (A2) as an
{\em NC-graph} (abbreviating ``normal crystal graph'').
  \Xcomment{
, using the term ``crystal'' for the pre-crystals
satisfying two additional (A3),(A4) (indicated later).
  }

The quadruple $(\ellin_1(v),\ellout_1(v),\ellin_2(v),\ellout_2(v))$
giving an important information about a vertex $v$ is called the
{\em length-tuple} of $v$ and denoted by $\tau(v)$.
We can also associate with a vertex $v$ the pair of integers
$\sigma(v):=(\ellout_1(v)-\ellin_1(v),\ellout_2(v)-\ellin_2(v))$.
From axiom (A2) it follows that
the difference $\sigma(v)-\sigma(u)$ is equal to $(-2,1)$ for all
1-edges $(u,v)$. In its turn, for
each 2-edge $(u,v)$, such a difference is equal to $(1,-2)$. So,
under the map $\sigma:V\to\Rset^2$, each 1-edge (2-edge) $e$ becomes
a parallel translation of the same vector $(-2,1)$ (resp. $(1,-2)$).
This implies the following property (indicated
for $n$-colored crystals in~\cite{ste-03}).

\begin{corollary} \label{cor:grad}
An NC-graph $K$ is graded w.r.t. each color $i$, in
the sense that any two directed paths in $K$ having the same first
vertex and the same last vertex contain equal numbers of $i$-edges.
\end{corollary}

In particular, $K$ is acyclic and has no parallel edges.
Also one can observe that the images by $\sigma$ of the monochromatic
subgraphs $(V,E_i)$ of $K$ have certain symmetries (assuming that
$\sigma$ is properly extended to the edges of $K$). If a 1-line has
length $p$, then $\sigma$ brings it to the directed straightline
segment in $\Rset^2$ with the beginning point $(p,q)$ and
the end point $(-p,p+q)$ for some $q$. Therefore, the median points
of the images of all 1-lines lie on the ``vertical'' coordinate axis
in $\Rset^2$. Similarly, the median points of the images of all
2-lines lie on the ``horizontal'' coordinate axis.

\medskip
\noindent{\bf Remark 1.}
NC-graphs have a rather loose structure, in contrast to RC-graphs.
In particular, the local finiteness (in the sense that all
monochromatic paths are finite, as required in axiom (A1)) and
even the finiteness of the list of different length-tuples of
vertices does not guarantee that the set of vertices is finite.
(We shall see later, in Remark~3 in Section~\ref{sec:theo}, that
a similar behavior is possible even if axiom (A3) is added.)
Indeed, for an arbitrary NC-graph $K$ having an undirected cycle,
one can construct an NC-graph with the same list of
length-tuples, which is viewed as an infinite tree
(the ``free covering'' over $K$). Also one can combine
NC-graphs as follows. Suppose a vertex $v$ of an NC-graph $K_1$
and a vertex $v'$ of an NC-graph $K_2$ have equal
lenght-tuples ($K_2$ may be taken as a copy of $K_1$).
Choose an edge in $K_1$ incident with $v$ and the corresponding edge
in $K_2$ incident with $v'$; let for definiteness these edges be
incoming 1-edges $(u,v)$ and $(u',v')$. Then the digraph obtained by
replacing these edges by $(u,v')$ and $(u',v)$ is also an NC-graph
(provided it is connected).

\medskip
To formulate two remaining axioms defining RC-graphs, we need some
definitions and notation.
(These axioms give combinatorial analogs of the Serre relations.)

In an NC-graph $K$, the edges with color $i$ ($i=1,2$) are naturally
associated with operator $F_i$ acting on the corresponding
subset of vertices. So, for a 1-edge (2-edge) $(u,v)$, we write
$v=F_1(u)$ and $u=F_1^{-1}(v)$ (resp. $v=F_2(u)$ and $u=F_2^{-1}(v)$).
Using this notation, one can express any vertex via another one (since
$K$ is connected). For example, the expression $F_1^{-1}F_2^2F_1(v)$
determines the vertex $w$ obtained from a vertex $v$ by traversing
1-edge $(v,v')$, followed by traversing 2-edges $(v',u)$
and $(u,u')$, followed by traversing 1-edge $(w,u')$ in backward
direction. Emphasize that every time we use an expression with $F_1$
or $F_2$ in what follows, this automatically says that all
involved edges do exist in $K$.

For each edge $e=(u,v)$ with color $i$, we assign {\em label}
$\ell(e):=0$ if $\ellout_{3-i}(u)=\ellout_{3-i}(v)$, and label
$\ell(e):=1$ otherwise. Axiom (A2) shows that the labels are
monotonically nondecreasing along any $i$-line $P$. In terms of
labels, the critical vertex in $P$ is just the vertex where the
incoming $i$-edge, if exists, is labeled 0 and the outgoing $i$-edge,
if exists, is labeled 1.

In further illustrations we will draw 1-edges by horizontal arrows
directed to the right, and 2-edges by vertical arrows directed up.

The first additional axiom describes situations when the
operators $F_1$ and $F_2$ commute.
\begin{itemize}
\item[(A3)] (a) If a vertex $u$ has two outgoing edges $(u,v),(u,v')$
and if $\ell(u,v)=0$, then $\ell(u,v')=1$ and $F_2F_1(u)=F_1F_2(u)$.
Symmetrically: (b) if a vertex $v$ has two incoming edges
$(u,v),(u',v)$ and if $\ell(u,v)=1$, then $\ell(u',v)=0$ and
$F_2^{-1}F_1^{-1}(v)=F_1^{-1}F_2^{-1}(v)$.
  \end{itemize}

Let us say that vertices $\tilde u,\tilde v,\tilde u',\tilde v'$
{\em form a square} if, up to renaming them,
$\tilde v=F_1(\tilde u)$, $\tilde u'=F_2(\tilde u)$ and $\tilde v'=
F_1(\tilde u')=F_2(\tilde v)$. The opposite 1-edges $(\tilde u,
\tilde v)$ and $(\tilde u',\tilde v')$ for this square have
equal labels, because of the obvious relations
$\ellin_2(\tilde u')=\ellin_2(\tilde u)+1$ and $\ellin_2(\tilde v')
=\ellin_2(\tilde v)+1$, and similarly for the opposite 2-edges
$(\tilde u,\tilde u')$ and $(\tilde v,\tilde v')$. Therefore,
(a) in (A3) implies $\ell(v,w)=1$ and $\ell(v',w)=0$,
where $w:=F_2F_1(u)$, and (b) implies $\ell(w,u)=0$ and
$\ell(w,u')=1$, where $w:=F_2^{-1}F_1^{-1}(v)$. The picture
illustrates the cases when edge $(u,v)$ has color 1.

 \begin{center}
  \unitlength=1mm
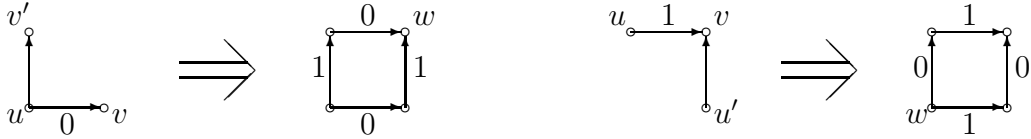
  \begin{picture}(140,20)
\put(5,5){\circle{1.0}}
\put(15,5){\circle{1.0}}
\put(45,5){\circle{1.0}}
\put(55,5){\circle{1.0}}
\put(95,5){\circle{1.0}}
\put(125,5){\circle{1.0}}
\put(135,5){\circle{1.0}}
\put(5,15){\circle{1.0}}
\put(45,15){\circle{1.0}}
\put(55,15){\circle{1.0}}
\put(85,15){\circle{1.0}}
\put(95,15){\circle{1.0}}
\put(125,15){\circle{1.0}}
\put(135,15){\circle{1.0}}
\put(5,5){\vector(1,0){9.5}}
\put(45,5){\vector(1,0){9.5}}
\put(125,5){\vector(1,0){9.5}}
\put(45,15){\vector(1,0){9.5}}
\put(85,15){\vector(1,0){9.5}}
\put(125,15){\vector(1,0){9.5}}
\put(5,5){\vector(0,1){9.5}}
\put(45,5){\vector(0,1){9.5}}
\put(55,5){\vector(0,1){9.5}}
\put(95,5){\vector(0,1){9.5}}
\put(125,5){\vector(0,1){9.5}}
\put(135,5){\vector(0,1){9.5}}
\put(25,9){\line(1,0){9}}
\put(25,11){\line(1,0){9}}
\put(105,9){\line(1,0){9}}
\put(105,11){\line(1,0){9}}
\put(31,6){\line(1,1){4}}
\put(31,14){\line(1,-1){4}}
\put(111,6){\line(1,1){4}}
\put(111,14){\line(1,-1){4}}
\put(2,3){$u$}
\put(2,16){$v'$}
\put(16,3){$v$}
\put(56,16){$w$}
\put(82,16){$u$}
\put(96,16){$v$}
\put(96,3){$u'$}
\put(121.5,3){$w$}
\put(9,1.5){0}
\put(49,1.5){0}
\put(129,1.5){1}
\put(49,16){0}
\put(89,16){1}
\put(129,16){1}
\put(42.5,9){1}
\put(56,9){1}
\put(122.5,9){0}
\put(136,9){0}
  \end{picture}
 \end{center}

From (A3) it follows that
\begin{numitem}
if $v$ is the critical vertex in an $i$-line, then $v$ is
simultaneously the critical vertex in the $(3-i)$-line passing
through $v$.
 \label{eq:crit}
\end{numitem}

Indeed, let for definiteness $i=1$ and assume that $v$ has outgoing
2-edge $(v,w)$. Suppose this edge is labeled 0. Then
$\ellin_1(v)>\ellin_1(w)$ implies that $v$ has incoming 1-edge
$(u,v)$. It is labeled 0 (since $v$ is the critical vertex in
$P_1(v)$). This means that $\ellout_2(u)=\ellout_2(v)>0$, and
therefore, $u$ has outgoing 2-edge $(u,u')$. Axiom (A3) implies
$w=F_1(u')$ and $\ell(u,u')=1$. But the latter contradicts the fact
that $(v,w)$ is labeled 0. Thus, $\ell(v,w)=1$. Arguing similarly and
using (b) in (A3), one shows that if $v$ has incoming 2-edge, then
this edge is labeled 0.

Therefore, we can speak of critical vertices without indicating line
colors. The final axiom indicates situations of ``remote
commuting'' $F_1$ and $F_2$.
\begin{itemize}
\item[(A4)] (i) If a vertex $u$ has two outgoing edges both
labeled 1, then $F_1F_2^2F_1(u)=F_2F_1^2F_2(u)$. Symmetrically:
(ii) if $v$ has two incoming edges both labeled 0, then
$F_1^{-1}(F_2^{-1})^2 F_1^{-1}(v)=F_2^{-1}(F_1^{-1})^2 F_2^{-1}(v)$.
  \end{itemize}

Note that in case (i), we have $F_2F_1(u)\ne F_1F_2(u)$ (otherwise
the vertices $u$, $v:=F_1(u)$, $v':=F_2(u)$ and $w:=F_2(v)$ would
form a square; then both edges $(v,w)$ and $(v',w)$ have label 1,
contrary to (A3)(b)). Similarly, $F_2^{-1}F_1^{-1}(v)\ne
F_1^{-1}F_2^{-1}(v)$ in case (ii).
The picture below illustrates axiom (A4). Here also the labels for
all involved edges are indicated and the critical vertices are
surrounded by circles.
(These labels and critical vertices are determined uniquely, which
is not difficult to show by use of (A3) and (A4). These facts will be
seen from the analysis in the next section as well.)

 \begin{center}
  \unitlength=1mm
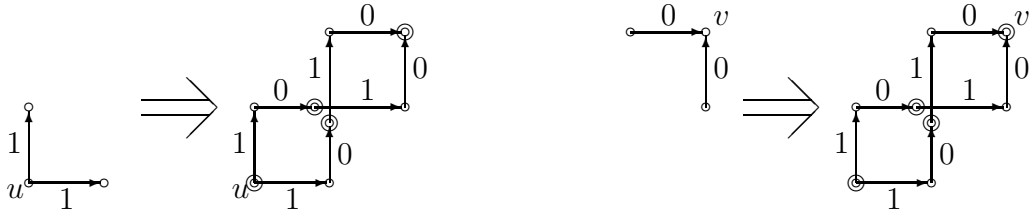
  \begin{picture}(147,30)
\put(5,5){\circle{1.0}}
\put(15,5){\circle{1.0}}
\put(5,15){\circle{1.0}}
\put(2,3){$u$}
\put(5,5){\vector(1,0){9.5}}
\put(5,5){\vector(0,1){9.5}}
\put(9,1.5){1}
\put(2,9){1}
\put(20,14){\line(1,0){9}}
\put(20,16){\line(1,0){9}}
\put(26,11){\line(1,1){4}}
\put(26,19){\line(1,-1){4}}
\put(35,5){\circle{1.0}}
\put(45,5){\circle{1.0}}
\put(35,15){\circle{1.0}}
\put(45,13){\circle{1.0}}
\put(43,15){\circle{1.0}}
\put(55,15){\circle{1.0}}
\put(45,25){\circle{1.0}}
\put(55,25){\circle{1.0}}
\put(35,5){\circle{2.0}}
\put(45,13){\circle{2.0}}
\put(43,15){\circle{2.0}}
\put(55,25){\circle{2.0}}
\put(35,5){\vector(1,0){9.5}}
\put(35,15){\vector(1,0){7.5}}
\put(43,15){\vector(1,0){11.5}}
\put(45,25){\vector(1,0){9.5}}
\put(35,5){\vector(0,1){9.5}}
\put(45,5){\vector(0,1){7.5}}
\put(45,13){\vector(0,1){11.5}}
\put(55,15){\vector(0,1){9.5}}
\put(32,3){$u$}
\put(39,1.5){1}
\put(32,9){1}
\put(46,7.5){0}
\put(37.5,16){0}
\put(49,16){1}
\put(42,19){1}
\put(56,19){0}
\put(49,26){0}
\put(85,25){\circle{1.0}}
\put(95,15){\circle{1.0}}
\put(95,25){\circle{1.0}}
\put(85,25){\vector(1,0){9.5}}
\put(95,15){\vector(0,1){9.5}}
\put(96,26){$v$}
\put(89,26){0}
\put(96,19){0}
\put(100,14){\line(1,0){9}}
\put(100,16){\line(1,0){9}}
\put(106,11){\line(1,1){4}}
\put(106,19){\line(1,-1){4}}
\put(115,5){\circle{1.0}}
\put(125,5){\circle{1.0}}
\put(115,15){\circle{1.0}}
\put(125,13){\circle{1.0}}
\put(123,15){\circle{1.0}}
\put(135,15){\circle{1.0}}
\put(125,25){\circle{1.0}}
\put(135,25){\circle{1.0}}
\put(115,5){\circle{2.0}}
\put(125,13){\circle{2.0}}
\put(123,15){\circle{2.0}}
\put(135,25){\circle{2.0}}
\put(115,5){\vector(1,0){9.5}}
\put(115,15){\vector(1,0){7.5}}
\put(123,15){\vector(1,0){11.5}}
\put(125,25){\vector(1,0){9.5}}
\put(115,5){\vector(0,1){9.5}}
\put(125,5){\vector(0,1){7.5}}
\put(125,13){\vector(0,1){11.5}}
\put(135,15){\vector(0,1){9.5}}
\put(136,26){$v$}
\put(119,1.5){1}
\put(112,9){1}
\put(126,7.5){0}
\put(117.5,16){0}
\put(129,16){1}
\put(122,19){1}
\put(136,19){0}
\put(129,26){0}
  \end{picture}
 \end{center}

Clearly the digraph obtained by reversing the orientation
of all edges of $K$, while preserving their colors, again satisfies
axioms (A1)--(A4) (thereby the label of each edge changes).
The resulting RC-graph is called {\em dual} to $K$ and is denoted
by $K^\ast$.

\section{Structural theorem}  \label{sec:theo}

In this section we present a theorem that clarifies the
combinatorial structure of RC-graphs defined by axioms (A1)--(A4).
According to this theorem, each RC-graph can be produced from two
elementary RC-graphs by use of a certain operation of replicating and
gluing together. First of all we introduce this operation in a general
form.

Consider arbitrary graphs or digraphs $G=(V,E)$ and $H=(V',E')$.
Let $S$ be a
distinguished subset of vertices of $G$, and $T$ a distinguished
subset of vertices of $H$.
Take $|T|$ disjoint copies of $G$, denoted as $G_t$ ($t\in T$), and
$|S|$ disjoint copies of $H$, denoted as $H_s$ ($s\in S$). We glue
these copies together in the following way: for each $s\in S$ and
each $t\in T$, the vertex $s$ in $G_t$ is identified with the vertex
$t$ in $H_s$. The resulting graph consisting of $|V| |T|+|V'| |S|-
|S| |T|$ vertices and $|E| |T|+|E'| |S|$ edges is denoted by
$(G,S)$\sprod$(H,T)$.

In our case the role of $G$ and $H$ play 2-colored digraphs
$K(a,0)$ and $K(0,b)$ depending on parameters $a,b\in\Zset_+$, each
of which being a certain triangular part of the Cartesian product of
two paths. More precisely, the vertices of $K(a,0)$ correspond to
the pairs $(i,j)$ for $i,j\in\Zset$ with $0\le i\le j\le a$, and the
vertices of $K(0,b)$ correspond to the pairs $(i,j)$ for
$0\le j\le i\le b$. The edges with color 1 in these graphs correspond to
all possible pairs of the form $((i,j),(i+1,j))$, and the edges with
color 2 to the pairs of the form $((i,j),(i,j+1))$. It is easy to check
that $K(a,0)$ satisfies axioms (A1)--(A4) and that the {\em diagonal}
$\{(i,i): i=0,\ldots,a\}$ is exactly the set of critical vertices in
it. Similarly, $K(0,b)$ is an RC-graph in which the set of critical
vertices coincides with the diagonal $\{(i,i): i=0,\ldots,b\}$.
These diagonals are just considered as the distinguished subsets $S$
and $T$ in these digraphs.

We refer to the digraph formed by applying the operation \sprod in
this case as the {\em diagonal-product} of $K(a,0)$ and $K(0,b)$, and
for brevity, denote it by $K(a,0)$\sprod$K(0,b)$. This
digraph is 2-colored, where the edge colors are inherited from
$K(a,0)$ and $K(0,b)$ in a natural way. The case $a=2$ and $b=1$ is
shown in the picture; here the critical vertices are marked with
circles. The trivial (degenerate) RC-graph $K(0,0)$ consists of a
unique vertex; clearly $K(a,0)$\sprod$K(0,0)=K(a,0)$ and
$K(0,0)$\sprod$K(0,b)=K(0,b)$ for any $a,b$.

\begin{figure}[htb]                
 \begin{center}
  \unitlength=1mm
  \begin{picture}(125,35)
\put(5,5){\circle{1.0}}           
\put(5,15){\circle{1.0}}
\put(5,25){\circle{1.0}}
\put(15,15){\circle{1.0}}
\put(15,25){\circle{1.0}}
\put(25,25){\circle{1.0}}
\put(5,5){\circle{2.5}}
\put(15,15){\circle{2.5}}
\put(25,25){\circle{2.5}}
\put(5,15){\vector(1,0){9.5}}
\put(5,25){\vector(1,0){9.5}}
\put(15,25){\vector(1,0){9.5}}
\put(5,5){\vector(0,1){9.5}}
\put(5,15){\vector(0,1){9.5}}
\put(15,15){\vector(0,1){9.5}}
\put(20,5){(a)}
\put(55,5){\circle{1.0}}
\put(65,5){\circle{1.0}}
\put(65,15){\circle{1.0}}
\put(55,5){\circle{2.5}}
\put(65,15){\circle{2.5}}
\put(55,5){\vector(1,0){9.5}}
\put(65,5){\vector(0,1){9.5}}
\put(70,5){(b)}
\put(95,5){\circle{1.0}}
\put(95,15){\circle{1.0}}
\put(95,25){\circle{1.0}}
\put(105,15){\circle{1.0}}
\put(105,25){\circle{1.0}}
\put(115,25){\circle{1.0}}
\put(95,5){\circle{2.5}}
\put(105,15){\circle{2.5}}
\put(115,25){\circle{2.5}}
\put(95,15){\vector(1,0){9.5}}
\put(95,25){\vector(1,0){9.5}}
\put(105,25){\vector(1,0){9.5}}
\put(95,5){\vector(0,1){9.5}}
\put(95,15){\vector(0,1){9.5}}
\put(105,15){\vector(0,1){9.5}}
\put(102,12){\circle{1.0}}
\put(102,22){\circle{1.0}}
\put(102,32){\circle{1.0}}
\put(112,22){\circle{1.0}}
\put(112,32){\circle{1.0}}
\put(122,32){\circle{1.0}}
\put(102,12){\circle{2.5}}
\put(112,22){\circle{2.5}}
\put(122,32){\circle{2.5}}
\put(102,22){\vector(1,0){9.5}}
\put(102,32){\vector(1,0){9.5}}
\put(112,32){\vector(1,0){9.5}}
\put(102,12){\vector(0,1){9.5}}
\put(102,22){\vector(0,1){9.5}}
\put(112,22){\vector(0,1){9.5}}
\put(102,5){\circle{1.0}}
\put(95,5){\vector(1,0){6.5}}
\put(102,5){\vector(0,1){6.5}}
\put(112,15){\circle{1.0}}
\put(105,15){\vector(1,0){6.5}}
\put(112,15){\vector(0,1){6.5}}
\put(122,25){\circle{1.0}}
\put(115,25){\vector(1,0){6.5}}
\put(122,25){\vector(0,1){6.5}}
\put(115,5){(c)}
  \end{picture}
 \end{center}
 \caption{(a) $K(2,0)$,\;\; (b) $K(0,1)$,\;\;
(c) $K(2,0)\triangleright\hspace{-4pt}\triangleleft\, K(0,1)$.}
  \label{fig:sprod}
  \end{figure}

It will be convenient for us to refer to a subgraph of an RC-graph
$K$ isomorphic to $K(a,0)$ (respecting colors and labels of the
edges), including $K(a,0)$ itself, as a {\em left sail}
of size $a$. Symmetrically, a subgraph of $K$ isomorphic to $K(0,b)$ is
referred to as a {\em right sail} of size $b$. In a left or right
sail we specify, besides the diagonal, the {\em 1-side} (the largest
1-line) and the {\em 2-side} (the largest 2-line).

It is a relatively easy exercise to verify validity of axioms
(A1)--(A4) for $K(a,0)$\sprod$K(0,b)$ with any $a,b\in\Zset_+$, i.e.
such a digraph is always an RC-graph. Our main theorem asserts that the
converse also takes place.

\begin{theorem} \label{tm:main}
Every RC-graph $K$ is representable as $K(a,0)$\sprod$K(0,b)$ for
some $a,b\in\Zset_+$. In particular, $K$ is finite.
  \end{theorem}

\noindent{\bf Proof of the theorem.}
The proof falls into several claims.

\medskip
\noindent{\bf Claim 1.} (i) {\em For any edge $(u,v)$ with color
$i$ and label 0, there exists edge $(w,u)$ with color $(3-i)$
and this edge has label 1.
Symmetrically: {\rm (ii)} for any $i$-edge $(u,v)$
labeled 1, there exists $(3-i)$-edge $(v,w)$ and this edge has label
0.}

\medskip
\begin{proof}
(i) For an $i$-edge $(u,v)$ labeled 0, one has
$\ellin_{3-i}(u)>\ellin_{3-i}(v)$ (by axiom (A2)). Therefore, $u$
has incoming $(3-i)$-edge $(w,u)$. Suppose $\ell(w,u)=0$.
Then $\ellout_i(w)=\ellout_i(u)>1$. So $w$ has outgoing $i$-edge
$(w,w')$. By axiom (A3) applied to the pair $(w,u),(w,w')$, the
vertices $w,w',u,v$ form a square, and $\ell(w,w')=1$. But the edge
$(u,v)$ opposite to $(w,w')$ for this square is labeled 0. This
contradiction shows that $(w,u)$ must be labeled 1.

Part (ii) in this claim follows from part (i) applied to the dual
RC-graph $K^\ast$.
\end{proof}

\medskip
\noindent{\bf Claim 2.} (i) {\em Let $(u,v)$ be an $i$-edge
labeled 0 and let $\ellout_{3-i}(v)>0$. Then there exist
$(3-i)$-edges $(u,u'),(v,v')$ labeled 1 and $i$-edges
$(u',v'),(v',v'')$ labeled 0.
Symmetrically: {\rm (ii)} if $(u,v)$ is an $i$-edge labeled 1 and if
$\ellin_{3-i}(u)>0$, then there exist $(3-i)$-edges
$(u',u),(v',v)$ labeled 0 and $i$-edges $(u'',u'),(u',v')$
labeled 1.}

\medskip
\begin{proof}
(i) For an $i$-edge $(u,v)$ labeled 0, one has
$\ellout_{3-i}(u)=\ellout_{3-i}(v)$. Therefore, $u$ has outgoing
$(3-i)$-edge $(u,u')$. By axiom (A3) applied to edges
$(u,v),(u,u')$, the vertices $u,v,u'$ and $v':=F_2(v)$ form a square,
and $\ell(u,u')=1$. Then $\ell(u',v')=\ell(u,v)=0$ and
$\ell(v,v')=\ell(u,u')=1$. The existence of $i$-ребра $(v',v'')$
labeled 0 follows from part (ii) in Claim~1 applied to the edge
$(v,v')$.

The second part of the claim follows from the first one applied
to $K^\ast$.
\end{proof}

\medskip
The picture below illustrates Claims 1 and 2 for the cases when
$(u,v)$ is a 1-edge labeled 0 or a 2-edge labeled 1.

 \begin{center}
  \unitlength=1mm
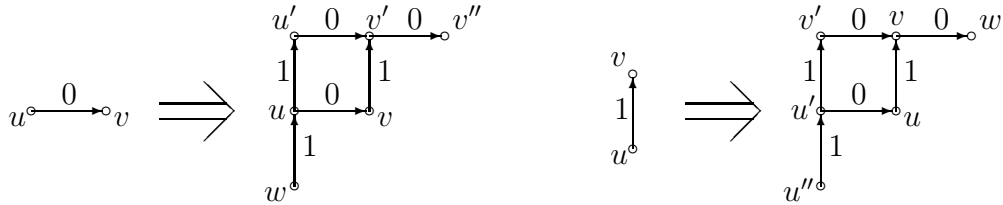
  \begin{picture}(140,30)
\put(5,15){\circle{1.0}}
\put(15,15){\circle{1.0}}
\put(5,15){\vector(1,0){9.5}}
\put(2,13){$u$}
\put(16,13){$v$}
\put(9,16){0}
\put(22,14){\line(1,0){9}}
\put(22,16){\line(1,0){9}}
\put(28,11){\line(1,1){4}}
\put(28,19){\line(1,-1){4}}
\put(40,5){\circle{1.0}}
\put(40,15){\circle{1.0}}
\put(40,25){\circle{1.0}}
\put(50,15){\circle{1.0}}
\put(50,25){\circle{1.0}}
\put(60,25){\circle{1.0}}
\put(40,15){\vector(1,0){9.5}}
\put(40,25){\vector(1,0){9.5}}
\put(50,25){\vector(1,0){9.5}}
\put(40,5){\vector(0,1){9.5}}
\put(40,15){\vector(0,1){9.5}}
\put(50,15){\vector(0,1){9.5}}
\put(36,3){$w$}
\put(36.5,14){$u$}
\put(37,26){$u'$}
\put(51,13){$v$}
\put(49.5,26){$v'$}
\put(61,26){$v''$}
\put(41,9){1}
\put(37.5,19){1}
\put(44,16){0}
\put(44,26){0}
\put(51,19){1}
\put(55,26){0}
\put(85,10){\circle{1.0}}
\put(85,20){\circle{1.0}}
\put(85,10){\vector(0,1){9.5}}
\put(82,8){$u$}
\put(82,21){$v$}
\put(82.5,14){1}
\put(92,14){\line(1,0){9}}
\put(92,16){\line(1,0){9}}
\put(98,11){\line(1,1){4}}
\put(98,19){\line(1,-1){4}}
\put(110,5){\circle{1.0}}
\put(110,15){\circle{1.0}}
\put(110,25){\circle{1.0}}
\put(120,15){\circle{1.0}}
\put(120,25){\circle{1.0}}
\put(130,25){\circle{1.0}}
\put(110,15){\vector(1,0){9.5}}
\put(110,25){\vector(1,0){9.5}}
\put(120,25){\vector(1,0){9.5}}
\put(110,5){\vector(0,1){9.5}}
\put(110,15){\vector(0,1){9.5}}
\put(120,15){\vector(0,1){9.5}}
\put(105,3){$u''$}
\put(106,14){$u'$}
\put(107,26){$v'$}
\put(121,13){$u$}
\put(119,26.5){$v$}
\put(131,26){$w$}
\put(111,9){1}
\put(107.5,19){1}
\put(114,16){0}
\put(114,26){0}
\put(121,19){1}
\put(125,26){0}
  \end{picture}
 \end{center}

For a path $(v_0,e_0,v_1,\ldots,e_k,v_k)$, we may use the abbreviate
notation $v_0v_1\ldots v_k$.

\medskip
\noindent{\bf Claim 3.} {\em Let $v$ be a critical vertex in $K$ and
let $L$ be a left sail of maximum size that contains $v$.
Then $L$ has size $d:=\ellin_1(v)+\ellout_2(v)$ and contains the paths
$\Pin_1(w)$ and $\Pout_2(w)$ for all vertices $w$ in $L$ (which is
equivalent to saying that $\ellout_2(w')=0$ for each vertex $w'$
on the 1-side of $L$, and $\ellin_1(w'')=0$ for each vertex
$w''$ on the 2-side of $L$). }

\smallskip
($L$ exists since the vertex $v$ itself forms the trivial sail
$K(0,0)$.)

\medskip
\begin{proof}
The claim is obvious if $d=0$. Let $d>0$. If $v$ has incoming 1-edge
$(u,v)$, then $\ell(u,v)=0$ (since $v$ is critical). By Claim~1, $v$
belongs to the left sail of size 1 formed by the edge $(u,v)$ and the
2-edge incoming $u$. Similarly, if $v$ has outgoing 2-edge, then it
belongs to a sail of size 1.

Thus, one may assume that the maximum-size left sail $L$ has size
$k\ge 1$. Consider the 1-side $P=v_0v_1\ldots v_k$ of $L$.
All 1-edges of $L$ are labeled 0, therefore,
$\ellout_2(v_0)=\ellout_2(v_1)=\ldots=\ellout_2(v_k)$.
Suppose $\ellout_2(v_i)>0$. Applying Claim~2
to the edges of $P$, one can conclude that there exists a path
$u_0u_1\ldots u_{k+1}$ whose edges have color 1 and label 0 and
whose vertices are connected with the vertices of $P$ by the 2-edges
$(v_i,u_i)$ labeled 1, $i=0,\ldots,k$. But this
implies that the sail $L$ is not maximum. Hence
$\ellout_2(v_i)=0$ для всех вершин $v_i$ в $P$.
Considering the 2-side of $L$ and arguing in a similar fashion, we
obtain $\ellin_1(w')=0$ for all vertices on this side, and the claim
follows.
  \end{proof}

Note that, in Claim~3, the vertex $v$ lies on the diagonal
of the sail $L$ (since the edges in $\Pout_1(v)$ are labeled 1).
One can see that $v$ determines $L$ uniquely, and therefore, we may
call $L$ the
{\em maximal left sail} containing $v$. Also each vertex $q$
in the diagonal of $L$ is critical. Indeed, suppose this is not so
for some $q$ and consider the 1-line $P$ passing through $q$ and the
critical vertex $r$ in $P$. Then $q$ occurs in $P$ earlier than $r$,
i.e., $q\in \Pin_1(r)$ and $q\ne r$. Applying Claim~3 to the
maximum-size left sail $L'$ containing $r$, one can conclude that
$L'$ includes $L$, contrary to the maximality of $L$.

Let $\Lscr$ be the set of all maximal left sails
(each containing a critical vertex). By above reasonings, the members
of $\Lscr$ are pairwise disjoint and they cover all
critical vertices, all 1-edges labeled 0 and all 2-edges labeled 1.

Since any left sail of $K$ turns into a right sail of $K^\ast$,
we have similar properties for the set $\Rscr$ of all
maximal right sails of $K$: the diagonal of each member of $\Rscr$
consists of critical vertices, the members of $\Rscr$ are pairwise
disjoint and they cover all critical vertices, all 1-edges
labeled 1 and all 2-edges labeled 0.

For a sail $Q$, the vertices in the diagonal $D$ of $Q$
are ordered in a natural way, where the minimal (maximal) element is
the vertex with zero indegree (resp. zero outdegree) in $Q$.
According to this ordering, the elements of $D$ are numbered by
$0,1,\ldots,|D|$. So if $u,v$ are vertices in $D$ with numbers
$i,i+1$, respectively, then $v=F_1F_2(u)$ in case of left sail, and
$v=F_2F_1(u)$ in case of right sail.

Choose a maximal left sail $L\in\Lscr$. Let it have size $a$ and
diagonal $D=(v_0,v_1,\ldots,v_a)$, where $v_i$ is the vertex with
number $i$ in $D$. For $i=0,\ldots,a$, denote by $R_i$ the maximal
right sail containing $v_i$. It has size $\ellout_1(v_i)+\ellin_2(v_i)$.
Since $v_{i+1}=F_1F_2(v_i)$ and since the 2-edge $(v_i,F_2(v_i))$ of
$L$ is labeled 1, we have $\ellout_1(v_{i+1})=\ellout_1(F_2(v_i))-1=
\ellout_1(v_i)$. In a similar way, the fact that the 1-edge
$(F_2(v_i),v_{i+1})$ of $L$ is labeled 0 implies
$\ellin_2(v_{i+1})=\ellin_2(v_i)$. Therefore, the values
$\ellout_1(v_i)$ ($i=0,\ldots,a$) are equal to one and the same number
$p$, and the values $\ellin_2(v_i)$ are equal to the same number $q$.
This gives the following important property:
  \begin{numitem}
the maximal right sails $R_0,\ldots,R_a$ have the same size $b:=p+q$,
and for $i=0,\ldots,a$, the vertex $v_i$ has number $q$ in the
diagonal of $R_i$ (and therefore, this number does
not depend on $i$).
  \label{eq:Ri}
  \end{numitem}
Then the sails $R_0,\ldots,R_a$ are
different, or, equivalently, $L$ and $R_i$ have a unique vertex in
common, namely, $v_i$ (since $R_i=R_j$ for $i\ne j$ would imply
$\ellin_2(v_i)\ne \ellin_2(v_j)$.)

Considering an arbitrary $R\in\Rscr$ and arguing similarly, we have:
  \begin{numitem}
the maximal left sails $L'\in\Lscr$
intersecting $R$ have the same size, and in the diagonals of these
sails $L'$, the vertices common with $R$ have equal numbers.
  \label{eq:L'}
  \end{numitem}

Since $K$ is connected, one can conclude from~\refeq{Ri}
and~\refeq{L'} that all members of $\Lscr$ have size $a$ and all
members of $\Rscr$ have size $b$.

Now let $u_0,\ldots,u_b$ be the vertices in the diagonal of the
maximal right sail $R_0$ (where the vertices are indexed according to
the ordering in the diagonal), and for $j=0,\ldots,b$, define $L_j$
to be the maximal left sail containing $u_j$. One may
assume that $L_0$ is just the left sail $L$ chosen above, i.e.,
$u_0=v_0$.

So far we have applied only axioms (A1)--(A3). The final
claim essentially uses axiom (A4).

\medskip
\noindent{\bf Claim 4.}
{\em For each $i=0,\ldots,a$ and each $j=0,\ldots,b$, the sails
$R_i$ and $L_j$ are intersecting. Moreover, their (unique) common
vertex $w$ has number $i$ in the diagonal of $L_j$ and has number
$j$ in the diagonal of $R_i$.
}

\medskip
\begin{proof}
We use induction on $i+j$. In fact, we have seen above that the claim
is valid for $j=0$ and any $i$, as well as for $i=0$ and any $j$.
So let $0<i\le a$ and $0<j\le b$. By induction there exist a common
vertex $u$ for $L_{j-1}$ and $R_{i-1}$, a common vertex $v$ for
$L_{j-1}$ and $R_i$, and a common vertex $v'$ for $R_{i-1}$ and $L_j$.
Furthermore, $u$ and $v$ have numbers $i-1$ and $i$ (respectively) in
the diagonal of $L_{j-1}$, and therefore, $v=F_2F_1(u)$. In their
turn, $u$ and $v'$ have numbers $j-1$ and $j$ (respectively) in the
diagonal of $R_{i-1}$, and therefore, $v'=F_2F_1(u)$. Also $v$ has
number $j-1<b$ in the diagonal of $R_i$, while $v'$ has number
$i-1<a$ in the diagonal of $L_j$. Hence the diagonal $D$ of $R_i$
contains vertex $w$ next to $v$ (i.e., $w=F_2F_1(v)$ and $w$ has
number $j$ in $D$), and the diagonal $D'$ of $L_j$ contains vertex
$w'$ next to $v'$ (i.e., $w'=F_1F_2(v')$ and $w'$ has number $i$ in
$D'$).

Since $(u,F_2(u))$ is a 2-edge of a left sail and $(u,F_1(u))$ is a
1-edge of a right sail, both edges are labeled 1. Applying axiom
(A4)(i) to them, we obtain $w=w'$, and the result follows.
\end{proof}

\medskip
Thus, the union $K'$ of sails $L_0,\ldots,L_b,R_0,\ldots,R_a$ is
isomorphic to $K(a,0)$\sprod$K(0,b)$. Also each of these sails meets
any other member of $\Lscr\cup\Rscr$ only within the set of critical
vertices in the former (by Claim~3), and these critical vertices
belong to $K'$. Therefore, the connectedness of $K$ implies $K=K'$.
This completes the proof of the theorem. $\qed$ $\qed$

\medskip
We denote the RC-graph isomorphic to $K(a,0)$\sprod$K(0,b)$ by
$K(a,b)$.

\medskip
  \begin{corollary} \label{cor:numbers}
$K(a,b)$ contains $(a+1)(b+1)$ critical vertices and has exactly
one vertex $s$ with zero indegree (the minimal vertex, or the
{\em source} of the RC-graph) and exactly one vertex $t$ with zero
outdegree (the maximal vertex, or the {\em sink}).
This $s$ is the common vertex of the sails $L_0$ and $R_0$ (defined
in the above proof), while $t$ is the common vertex of the sails
$L_b$ and $R_a$, and one holds $\ellout_1(s)=\ellin_2(t)=b$ and
$\ellout_2(s)=\ellin_1(t)=a$. In particular, an RC-graph with source
$s$ is determined by the parameters $\ellout_1(s),\ellout_2(s)$. Also
$K(a,b)$ has equal numbers of edges of each color, namely,
$(a+1)(b+1)(a+b)/2$.
  \end{corollary}

\noindent{\bf Remark 2.}
In the proof of Theorem~\ref{tm:main} we never used part (ii) of
axiom (A4). This means that this part is redundant. In fact, one can
show directly that (A4)(ii) is implied by (A4)(i) and (A1)--(A3).
(However, part (ii) of axiom (A4) becomes essential for infinite
RC-graphs considered in the next section.)

\medskip
\noindent{\bf Remark 3.}
The connected 2-colored digraphs satisfying axioms
(A1)--(A3), which may be named {\em weakened RC-graphs}, or
{\em WC-graphs},
form an interesting class $\Wscr\Cscr$ lying between the classes of
NC-graphs and RC-graphs.
This class can be completely characterized, relying on the
fact (seen from the above proof)
that each member of it is the union of a set of pairwise
disjoint left sails of the same size $a$ and a set of pairwise
disjoint right sails of the same size $b$. More precisely, each
$K\in\Wscr\Cscr$ can be encoded by
parameters $a,b\in\Zset_+$, a graph $\Gamma=(V(\Gamma),E(\Gamma))$
and a map $\omega:V(\Gamma)\to \Zset_+$ such that:
  \begin{itemize}
\item[($\ast$)] (a) $\Gamma$ is a connected finite
or infinite bipartite graph with vertex parts $V_1,V_2$;
(b) each vertex in $V_1$ has degree $a+1$ and each vertex in $V_2$
has degree $b+1$; (c) $\omega(v)$ does not exceed $b$ for each
vertex in $V_1$, and $a$ for each vertex in $V_2$; (d) for each
$v\in V(\Gamma)$, the vertices $w$ adjacent to $v$ have
different values $\omega(w)$.
  \end{itemize}
The vertices in $V_1$ are associated
with left sails of size $a$, and the vertices in $V_2$ with right
sails of size $b$. The edges of $\Gamma$ indicate how these sails are
glued together, namely: if $u\in V_1$ and $v\in V_2$ are connected by an
edge, then the vertex with number $\omega(v)$ in the diagonal of the
sail corresponding to $u$ is identified with the vertex with number
$\omega(u)$ in the diagonal of the sail corresponding to $v$. The
resulting $K$ is finite if and only if $\Gamma$ is finite. If
$\Gamma$ is a complete bipartite graph, we just obtain the RC-graph
$K(a,b)$. A somewhat different way to characterize $\Wscr\Cscr$ is as
follows.

For $a,b\in\Zset_+$, let $\Gamma_{a,b}$ denote the
Cartesian product $P_a\times P_b$ of directed paths with length
$a$ and $b$, respectively.
A {\em covering} over $\Gamma_{a,b}$ is a nonempty (finite or
infinite) connected digraph $G$ along with a homomorphism
$\gamma:G\to\Gamma_{a,b}$ under which the 1-neighborhood of each
vertex $v$ of $G$ (i.e., the subgraph induced by the edges incident
with $v$) is isomorphically mapped to the 1-neighborhood of
$\gamma(v)$. For such a $(G,\gamma)$, the preimage in $G$ of each path
$(P_a,\cdot)$ of $\Gamma_{a,b}$ is a collection of pairwise disjoint
paths $Q$ of length $a$, and we can replace each of these $Q$
by a copy of the left sail $L$ of size $a$ in a natural way (the
vertices of $Q$ are identifyed with the elements of the diagonal of
$L$ in the natural order). The preimages of each path $(\cdot,P_b)$
are replaced by copies of the right sail $R$ of size $b$ in a similar
fashion. One can check that the resulting digraph $K$, with a due
assignment of edge colors, satisfy (A1)--(A3), and conversely, each
member of $\Wscr\Cscr$ can be obtained by this construction.

  \begin{prop} \label{pr:wc}
There is a bijection between the set $\Wscr\Cscr$ of 2-colored
digraphs satisfying axioms (A1)--(A3) and the set of
coverings over grids $\Gamma_{a,b}$ for all $a,b\in\Zset_+$.
  \end{prop}
Adding axiom (A4) removes all nontrivial coverings (i.e.,
those different from the grids themselves). The fact that the grid
$\Gamma_{a,b}$ has one source and one sink (the zero indegree and
zero outdegree vertices, respectively) implies that the quantity of
sources (sinks) in a covering $G$ is equal to the quantity of
preimages in $G$ of a vertex of $\Gamma_{a,b}$. This gives the
following important property.

  \begin{corollary} \label{cor:1source}
Under validity of (A1)--(A3), axiom (A4) is equivalent to the
requirement that the digraph has only one source or only one sink. In
other words, the RC-graphs are precisely the WC-graphs with one
source (one sink).
  \end{corollary}

\noindent{\bf Remark 4.}
The fact that an RC-graph $K=K(a,b)$ is graded w.r.t. each color
(cf. Corollary~\ref{cor:grad}) implies
that $K$ is the Hasse diagram of a poset $(V,\preceq)$
on the vertex set (it is generated by the relations $u\prec v$ for
edges $(u,v)$). Considering the sail structure of $K$ as above, it is not
difficult to obtain the following sharper property. Here $r_{ij}$
stands for the common critical vertex of sails $R_i$ and $L_j$, and
for a vertex $v$, we denote by $p(v)$ ($q(v)$) the minimal (resp.
maximal) critical vertex greater (resp. smaller) than or equal to $v$
in a maximal sail containing $v$.

  \begin{prop} \label{pr:lattice}
The poset $(V,\preceq)$ is a lattice, that is, any two vertices $u,v$
have a unique minimal upper bound $u\vee v$ and a unique maximal
lower bound $u\wedge v$. More precisely:

{\rm (i)} if $u,v$ are critical vertices $r_{ij},r_{i'j'}$, then
$u\vee v=r_{\max\{i,i'\},\max\{j,j'\}}$ and
$u\wedge v=r_{\min\{i,i'\},\min\{j,j'\}}$;

{\rm (ii)} if $u,v$ occur in the same maximal sail, then both
$u\vee v$ and $u\wedge v$ belong to this sail (they are
computed in a straightforward way; in particular, $p(u\vee v)=
p(u)\vee p(v)$ and $q(u\vee v)=q(u)\vee q(v)$);

{\rm (iii)} for vertices $u,v$ occurring in different maximal sails:
(a) if $p(u)\preceq p(v)$, then $v$ and $w:=u\vee v$ belong to the
same maximal sail (sails) $Q$ and one holds $p(w)=p(v)$ and
$q(w)=p(u)\vee q(v)$ (the latter vertex belongs to $Q$ as well);
and (b) if $p(u)$ and $p(v)$ are incomparable,
then $u\vee v=p(u)\vee p(v)$.
  \end{prop}
(For vertices $u,v$ in different maximal sails, computing
$u\wedge v$ is symmetric to (iii).) Note that this lattice is not
distributive already for $a,b=1$.

\section{Infinite RC-graphs} \label{sec:infin}

The notion of RC-graphs can be extended, with a due care, to
(connected) infinite 2-colored digraphs $K=(V,E_1,E_2)$.
By Theorem~\ref{tm:main}, the finiteness of each monochromatic
path leads to the finiteness of an RC-graph, so we now should
allow infinite monochromatic paths and accordingly modify axiom (A1).
There are three types of infinite paths. A {\em fully infinite path}
is a sequence of the form $\ldots,v_i,e_i,v_{i+1},e_{i+1},\ldots$,
where the index set $I$ of vertices $v_i$ ranges $\Zset$
(and, as before, $e_i$ is the edge from $v_i$ to $v_{i+1}$). If
$I=\Zset_+$ (resp. $I=\Zset_-$), we deal with a {\em semiinfinite
path} in forward (resp. backward) direction.

There are two methods to define an RC-graph so as to involve both
finite or infinite cases. The first method is based on a
generalization of Theorem~\ref{tm:main} and uses the construction
$(G,S)$\sprod$(H,T)$ defined in Section~\ref{sec:theo} and applicable
to arbitrary finite or infinite graphs $G,H$ and distinguished
subsets $S,T$.

In our case the role of $G$ plays a left sail, which can be either
finite (defined as before) or infinite of any of three possible
sorts. The sail denoted by $\Lscr^\infty_\infty$ (``infinite up
and to the left'') has vertex set $\{(i,j): \;i,j\in\Zset,\; i\le
j\}$, the sail $\Lscr^\infty$ (``infinite up'') has the vertices
$(i,j)$ for $0\le i\le j$, and the sail $\Lscr_\infty$
(``infinite to the left'') has the vertices $(i,j)$ for $i\le j\le 0$.
Analogously, $H$ is a finite of infinite right sail, and the latter
can be of three sorts: $\Rscr^\infty_\infty$ (``infinite down and to
the right'') with the vertices $(i,j)$ for $j\le i$,
$\Rscr^\infty$ (``infinite to the right'') with the vertices
$(i,j)$ for $0\le j\le i$, and $\Rscr_\infty$ (``infinite down'')
with the vertices $(i,j)$ for $j\le i\le 0$. In all cases the 1-edges
correspond to the pairs $((i,j),(i+1,j))$, and 2-edges to the pairs
$((i,j),(i,j+1))$. As before, in a left sail all 1-edges are labeled
0, and all 2-edges are labeled 1, while in a right sail the labels
are interchanged. The distinguished subsets $S,T$ are the
corresponding ``diagonals'' consisting of the vertices $(i,i)$
(being critical). Infinite sails are illustrated in
the picture.

 \begin{center}
  \unitlength=1mm
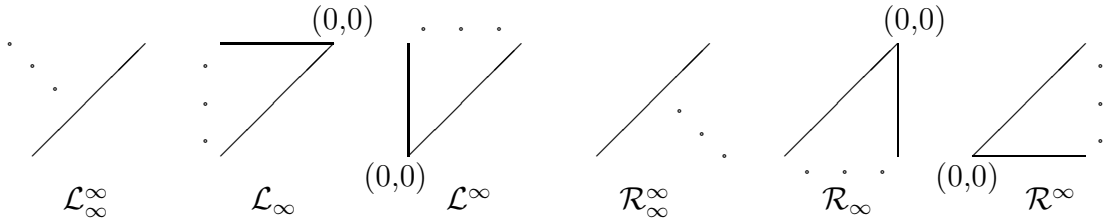
  \begin{picture}(150,30)
\put(5,10){\line(1,1){15}}
\put(2,25){\circle{0.5}}
\put(5,22){\circle{0.5}}
\put(8,19){\circle{0.5}}
\put(9,3){$\Lscr^\infty_\infty$}
\put(30,10){\line(1,1){15}}
\put(30,25){\line(1,0){15}}
\put(28,22){\circle{0.5}}
\put(28,17){\circle{0.5}}
\put(28,12){\circle{0.5}}
\put(34,3){$\Lscr_\infty$}
\put(42,27){(0,0)}
\put(55,10){\line(1,1){15}}
\put(55,10){\line(0,1){15}}
\put(57,27){\circle{0.5}}
\put(62,27){\circle{0.5}}
\put(67,27){\circle{0.5}}
\put(60,3){$\Lscr^\infty$}
\put(49,6){(0,0)}
\put(80,10){\line(1,1){15}}
\put(97,10){\circle{0.5}}
\put(94,13){\circle{0.5}}
\put(91,16){\circle{0.5}}
\put(83,3){$\Rscr^\infty_\infty$}
\put(105,10){\line(1,1){15}}
\put(120,10){\line(0,1){15}}
\put(108,8){\circle{0.5}}
\put(113,8){\circle{0.5}}
\put(118,8){\circle{0.5}}
\put(110,3){$\Rscr_\infty$}
\put(118,27){(0,0)}
\put(130,10){\line(1,1){15}}
\put(130,10){\line(1,0){15}}
\put(147,12){\circle{0.5}}
\put(147,17){\circle{0.5}}
\put(147,22){\circle{0.5}}
\put(137,3){$\Rscr^\infty$}
\put(125,6){(0,0)}
  \end{picture}
 \end{center}

Combining any of the four sorts of left sails (one of which is
finite and the other three are infinite)  with any of the four
sorts of right sails, we obtain 16 types of RC-graphs $L$\sprod$R$, of
which one is finite, while the other 15 contain a fully infinite or
semiinfinite monochromatic paths.
For example, $\Lscr_\infty$\sprod$\Rscr^\infty$ has fully infinite
1-lines and finite 2-lines, and $K(a,0)$\sprod$\Rscr^\infty_\infty$
has semiinfinite in forward direction 1-lines and semiinfinite in
backward direction 2-lines. The largest RC-graph
$\Lscr^\infty_\infty$\sprod$\Rscr^\infty_\infty$ contains the other
ones as induced subgraphs.

\medskip
\noindent
{\bf Remark 5.}
There are five more infinite RC-graphs (defined up to swapping the
edge colors). They have a simple structure and do not contain
critical vertices at all. The vertices of these
RC-graphs are the pairs $(i,j)$, where either (i) $i,j\in\Zset$, or
(ii) $i\in\Zset$ and $j\in\Zset_+$, or (iii) $i\in\Zset$ and
$j\in\Zset_-$, or (iv) $i\in\Zset_-$ and $j\in\Zset_+$, or
(v) $i\in\Zset_+$ and $j\in\Zset_-$. Formally, in cases (ii),(iv),
the 1-edges $((i,j),(i+1,j))$ are labeled 1, and 2-edges
$((i,j),(i,j+1))$ are labeled 0, while in cases (iii),(v), the labels
are interchanged.

\medskip
The second, alternative, method of unfying the definition of RC-graphs
to include infinite cases consists in modifying axioms (A2) and (A3)
(while preserving (A4)). We replace them by a single axiom that
postulates properties exposed in Claims~1 and~2 from the proof of
Theorem~\ref{tm:main}. As before, each edge $e$ is endowed with
label $\ell(e)\in\{0,1\}$, the labels are monotonically nondecreasing
along each monochromatic path, and $\ellin_i(v),\ellout_i(v)$ denote
the lengths of corresponding paths (which may be infinite).
The new axiom is stated as follows:
\begin{itemize}
\item[(A$'$)] $K$ is graded w.r.t. each color (cf.
Corollary~\ref{cor:grad}). Also: (a) for each $i$-edge $(u,v)$
labeled 0, there exists $(3-i)$-edge $(w,u)$ labeled 1; moreover,
$u$ has outgoing $(3-i)$-edge $(u,u')$ if and only if $v$ has
outgoing $(3-i)$-edge $(v,v')$, and in this case both
$(u,u'),(v,v')$ are labeled 1 and there exist $i$-edges
$(u',v'),(v',v'')$ labeled 0. Symmetrically: (b) for
each $i$-edge $(u,v)$ labeled 1, there exists $(3-i)$-edge $(v,w)$
labeled 0; moreover,
$u$ has incoming $(3-i)$-edge $(u',u)$ if and only if $v$ has
incoming $(3-i)$-edge $(v',v)$, and in this case both
$(u',u),(v',v)$ are labeled 0 and there exist $i$-edges
$(u',v'),(u'',u')$ labeled 1.
   \end{itemize}

In reality both definitions are equivalent.
A verification that a generalized RC-graph constructed by the
first method satisfies (A$'$) (and (A4)) is relatively simple.
The converse assertion, that a digraph
satisfying (A$'$),(A4) is one of those described in the first
method, can be proved by following the method of proof of
Theorem~\ref{tm:main}, with necessary extensions and refinements;
we omit details here. (Unlike the finite case (cf. Remark~2), part
(ii) of axiom (A2) becomes essential for the general case.)
In particular, one shows that if some (infinite) monochromatic
line has no critical vertex, then the RC-graph is one of those
indicated in Remark~5.

\medskip
\noindent
{\bf Remark 6.}
The construction of diagonal-product can be used for extending the
notion of RC-graphs to more abstract structures. More precisely,
let $I,J$ be two fully ordered sets. (For example, we can take as
$I,J$ intervals in $\Rset$ or $\Qset$. In essense, so far
we have dealt with intervals in $\Zset$.) We define the left sail
over $I$ in a natural way, to be the set $L:=\{(x,y)\in I^2:
x\preceq y\}$, and define the right sail over $J$ to be
$R:=\{(x,y)\in J^2: x\succeq y\}$. The distinguished subsets
$D,D'$, or the diagonals, in $L,R$, respectively, consist of the
identical pairs $(x,x)$. Then we can form
the corresponding ``diagonal-product'' $K:=L$\sprod$R$.
Fixing the second coordionate $y$ (resp. the first coordinate $x$)
in the left sail $L$ or in the right sail $R$ gives a line of
color 1 (resp. 2) in this sail. When $I,J$ are intervals in $\Rset$
(in which case $K$ may be named a {\em continuous crystal}),
one can introduce a reasonable metric on $K$, which determines
its intrinsic topological structure, as follows.
The distance between points $(x,y)$ and $(x',y')$ in each sail of $K$
is assigned to be the $\ell_1$-distance $|x-x'|+|y-y'|$, and the
distance within the critical set $W:=D\times D'$ is also assigned
to be the corresponding distance of $\ell_1$-type. This induces a
metric $d$ on the entire $K$:
for different sails $Q,Q'$ of $K$ and points $u\in Q$ and $v\in Q'$,
$d(u,v)$ is equal to $\inf\{d(u,u')+d(u',v')+d(v',v):
\;u'\in Q\cap W,\; v'\in Q'\cap W\}$. Note that the resulting metric
space need not be compact even if the intervals $I,J$ are
bounded and closed.

\section{Polyhedral aspects and a relation to Gelfand-Tsetlin
patterns} \label{sec:patt}

In this section we return to a (finite) RC-graph $K=K(a,b)$, with
node set $V$, and discuss some natural embeddings of $K$ and other
properties, using definitions, notation and results from
Section~\ref{sec:theo}. (The results can be extended to infinite
RC-graphs as well.)

Recall that $K$ has set $\Lscr=\{L_0,\ldots,L_b\}$ of maximal
left sails (with size $a$) and set $\Rscr=\{R_0,\ldots,R_a\}$ of
maximal right sails (with size $b$). These sails and the
vertices in their diagonals are numbered as in the proof of
Theorem~\ref{tm:main}. Under these numerations, sails
$L_j$ and $R_i$ intersect at the critical vertex that
has number $i$ in the diagonal $D_j$ of $L_j$ and number $j$ in
the diagonal $D'_i$ of $R_i$. We denote the vertex with
number 0 in $D_j$  (the {\em source} of $L_j$) by $s_j$, and denote
the vertex with number 0 in $D'_i$ (the source of $R_i$) by $s'_i$.
Each vertex $v$ of $L_j$ is determined by {\em (local) coordinates}
$(p,q)$, where $p$ (resp. $q$) is the number of 1-edges (resp.
2-edges) in a path from $s_j$ to $v$, i.e., $v=F_1^pF_2^q(s_j)$.
Analogous coordinates are assigned in $R_i$ with respect to
$s'_i$. The vertex $s_0=s'_0$ is the source of the whole $K$, denoted
by $s_K$.

\medskip
\noindent
{\bf 1.}
One way to embed $K$ in an Abelian group relyes on the observation
that the vertices $v\in V$ have different length-tuples $\tau(v)=
(\ellin_1(v),\ellout_1(v),\ellin_2(v),\ellout_2(v))$. Moreover, the
vertices differ from each other even if three parameters involved
in $\tau$ are considered, e.g.,
$\ellin_1(v),\ellout_1(v),\ellin_2(v)$. This is seen from the
following lemma.

\begin{lemma} \label{lm:eps}
For $v\in V$, define $\eps:=a-\ellin_1(v)-\ellout_2(v)$ and
$\delta:=b-\ellin_2(v)-\ellout_1(v)$.
{\rm (i)} If $v$ occurs in a left sail $L_j$, then $-2\delta=
\eps\ge 0$, $v$ has coordinates $(\ellin_1(v),a-\ellout_2(v))$ in
$L_j$, and $j$ is equal to $\ellin_2(v)-\eps=b-\ellout_1(v)+\eps$.
 {\rm (ii)} If $v$ occurs in a right sail $R_i$, then $-2\eps=
\delta\ge 0$, $v$ has coordinates $(b-\ellout_1(v),\ellin_2(v))$ in
$R_i$, and $i$ is equal to $\ellin_1(v)-\delta=a-\ellout_2(v)+\delta$.
(iii) The vertex $v$ is critical if and only if $\eps=\delta=0$.
  \end{lemma}
  \begin{proof}
The assertions are obvious when $v$ is critical. If $v$ lies in a
left sail $L_j$, then the assertions in (i) can be obtained by
comparing $\tau(v)$ with the length-tuples of the critical vertices
in the lines $P_1(v)$ and $P_2(v)$ and by using the fact that both
critical vertices have number $j$ in the diagonals of the
corresponding maximal right sails. If $v$ lies in a right sail $R_i$,
the proof is analogous.
  \end{proof}

Note that the edges $(u,v)$ with the same color and the same
label have the same difference $\tau(v)-\tau(u)$ (e.g., for color 1
and label 0, the difference is $(1,-1,-1,0)$). So $\tau$ induces an
embedding of $K$ in the corresponding subgroup of $\Zset^4$ shifted
by the vector $-\tau(s_K)=(0,-b,0,-a)$.

\medskip
\noindent
{\bf 2.}
Next we are interested in embeddings with the property that the
edges of $K$ correspond to parallel translations of unit base vectors.
For $l=0,1$ and a path $P$ in $K$, we denote the number of 1-edges
(2-edges) of $P$ with label $l$ by $\alpha_l(P)$ (resp. by
$\beta_l(P)$). The next lemma strengthens Corollary~\ref{cor:grad},
showing that $K$ is graded w.r.t. each combination of color and
label.

\begin{lemma}  \label{lm:2grad}
Let $P$ be a path in $K$ beginning at $s_K$ and ending at $v\in V$.
{\rm (i)} If $v$ occurs in a left sail $L_j$, then $\alpha_1(P)=
\beta_0(P)=j$, and $v$ has coordinates $(\alpha_0(P),\beta_1(P))$ in
$L_j$.
{\rm (ii)} If $v$ occurs in a right sail $R_i$, then $\alpha_0(P)=
\beta_1(P)=i$, and $v$ has coordinates $(\alpha_1(P),\beta_0(P))$ in
$R_i$.
  \end{lemma}
  \begin{proof}
Use induction on the length $|P|$ of $P$. The assertion is trivial
when $|P|=0$, so let $|P|>0$. Suppose $v$ lies in $L_j$. If
the vertex $u$ of $P$ preceding $v$ also lies in $L_j$, the assertion
for $P$ easily follows by induction from that for the part of $P$ from
$s_K$ to $u$.

Now let $u\not\in L_j$. Then the vertex $v$ is critical and the edge
$(u,v)$ is contained in some right sail $R_i$. Consider the last
critical vertex $w$ of $P$ different from $v$ (it exists as the
beginning vertex $s_K$ of $P$ is critical and $s_K\ne v$). Clearly
$w$ belongs to $R_i$; let it have number $j'$ in the diagonal of $R_i$
(whereas $v$ has number $j$, and $j>j'$). For the part $P'$ of $P$
from $s_K$ to $w$, one has $\alpha_0(P')=\alpha_0(P)$,
$\beta_1(P')=\beta_1(P)$, $\alpha_1(P')+j-j'=\alpha_1(P)$ and
$\beta_0(P')+j-j'=\beta_0(P)$. By induction $\alpha_1(P')=\beta_0(P')
=j'$ and $\alpha_0(P')=\beta_1(P')=i$ (since $w$ is also contained in
the left sail $L_{j'}$ and has number $i$ in its diagonal, whence $w$
has coordinales $(i,i)$ in $L_{j'}$). This gives the desired result
for $P$, taking into account that $v$ has the same coordinates
$(i,i)$ in $L_j$.

When $v$ lies in a right sail $R_i$, we argue in a similar way.
  \end{proof}

Thus, for a path $P$ from the source $s_K$ to a vertex $v$, the
numbers $\alpha_l(P),\beta_l(P)$ ($l=0,1$) depend only on $v$,
and we can define $\alpha_l(v):=\alpha_l(P)$ and
$\beta_l(v):=\beta_l(P)$. Also Lemma~\ref{lm:2grad} shows that
the quadruples $\bar\rho(v):=(\alpha_0(v),\alpha_1(v),\beta_0(v),
\beta_1(v))$ are different for all vertices $v$, i.e., the
map $\bar\rho:V\to\Zset^4$ is injective. Under this map, traversing
an edge of $K$ corresponds to adding a unit base vector associated
with the color and label of the edge.

Since the local coordinates $(p,q)$ satisfy the relation
$0\le p\le q\le a$ for the sails in $\Lscr$, and $0\le q\le p\le b$
for the sails in $\Rscr$, Lemma~\ref{lm:2grad} implies the following.

 \begin{corollary} \label{cor:neq}
For each vertex $v$, one has $0\le\alpha_0(v)\le\beta_1(v)\le a$ and
$0\le\beta_0(v)\le\alpha_1(v)\le b$; moreover, at least one of
$\alpha_0(v)\le\beta_1(v)$ and $\beta_0(v)\le\alpha_1(v)$ turns into
equality (and both equalities here characterize the critical
vertices). Conversely, if integers $p,p',q,q'$ satisfy
$0\le p\le q\le a$ and $0\le q'\le p'\le b$ and if at least one of
$p=q$ and $p'=q'$ holds, then there is a vertex $v$ with
$\bar\rho(v)=(p,p',q,q')$.
  \end{corollary}

\noindent {\bf Remark 7.}
Using Lemma~\ref{lm:2grad}, one can characterize the lattice
$(V,\preceq)$ of $K$ (cf. Proposition~\ref{pr:lattice}) via the
vertex parameters $\alpha_i,\beta_j$. More precisely, vertices $u,v$
satisfy $u\preceq v$ if and only if at least one of the following
holds:

{\rm (a)} $\beta_1(u)\le\alpha_0(v)$ and $\alpha_1(u)\le \beta_0(v)$;

{\rm (b)} all $\alpha_0(u),\alpha_0(v),\beta_1(u),\beta_1(v)$ are
equal, $\alpha_1(u)\le\alpha_1(v)$, and $\beta_0(u)\le\beta_0(v)$.

{\rm (c)} all $\alpha_1(u),\alpha_1(v),\beta_0(u),\beta_0(v)$ are
equal, $\alpha_0(u)\le\alpha_0(v)$, and $\beta_1(u)\le\beta_1(v)$.

\medskip
(One can describe the lattice operations $\vee,\wedge$ in terms of
$\alpha_i,\beta_j$; we leave this to the reader as an exercise.)

\bigskip
\noindent
{\bf 3.} Corollary~\ref{cor:neq} enables us to transform the map
$\bar\rho$ defined in part 2 into an injective map
$\rho:V\to\Zset^3$, by combining $\beta_0$ and $\beta_1$ into one
coordinate. More precisely, define $\beta:=\beta_0+\beta_1$ and
$\rho:=(\alpha_0,\alpha_1,\beta)$. The fact that $\rho$ is injective
follows from the possibility of (uniquely) restoring
$\beta_0(v),\beta_1(v)$ if we know $\rho(v)$, namely:
\begin{numitem}
$\beta_0(v)=\alpha_1(v)$ if $\beta(v)\ge \alpha_0(v)+\alpha_1(v)$,
and $\beta_1(v)=\alpha_0(v)$ otherwise; equivalently:
$\beta_0(v)=\min\{\alpha_1(v),\beta(v)-\alpha_0(v)\}$ and
$\beta_1(v)=\max\{\alpha_0(v),\beta(v)-\alpha_1(v)\}$.
  \label{eq:rest}
  \end{numitem}
\indent In a similar way, one can combine $\alpha_0$ and $\alpha_1$,
by setting $\alpha:=\alpha_0+\alpha_1$ and
$\rho':=(\alpha,\beta_0,\beta_1)$. Then the injectivity of $\rho'$
is provided by:
\begin{numitem}
$\alpha_0(v)=\beta_1(v)$ if $\alpha(v)\ge \beta_0(v)+\beta_1(v)$,
and $\alpha_1(v)=\beta_0(v)$ otherwise; equivalently:
$\alpha_0(v)=\min\{\beta_1(v),\alpha(v)-\beta_0(v)\}$ and
$\alpha_1(v)=\max\{\beta_0(v),\alpha(v)-\beta_1(v)\}$.
  \label{eq:rest1}
  \end{numitem}
\indent Consider the map $\rho$ and inentify $V$ with the set
$\rho(V)$ of points in the space $\Rset^3$ with coordinates
$(\alpha_0,\alpha_1,\beta)$. Let $\Pscr=\Pscr(a,b)$ denote the convex
hull of $V$. Using Corollary~\ref{cor:neq}, it is not difficult to
obtain the following description and properties of the polytope
$\Pscr$.

\begin{prop}  \label{pr:Pscr}
$P$ is formed by the vectors $(\alpha_0,\alpha_1,\beta)\in\Rset^3$
satisfying
\begin{eqnarray}
{\rm (i)} &\quad& 0\le \alpha_0\le a; \label{eq:poly}\\
{\rm (ii)} &\quad& 0\le \alpha_1\le b; \nonumber\\
{\rm (iii)} &\quad& \alpha_0\le \beta\le \alpha_1+a. \nonumber
  \end{eqnarray}
\noindent The polytope $\Pscr$ is represented as the Minkowsky sum of
the convex hulls of sails $L_0$ and $R_0$ (considered as sets of
points) and the set of integer points in $\Pscr$ is exactly $V$.
The vertices of $\Pscr$ are
$(0,0,0),(0,0,a),(0,b,0),(a,0,a),(a,b,a),(0,b,a+b),(a,b,a+b)$ (some
of which coincide when $a=0$ or $b=0$).
\qed
  \end{prop}

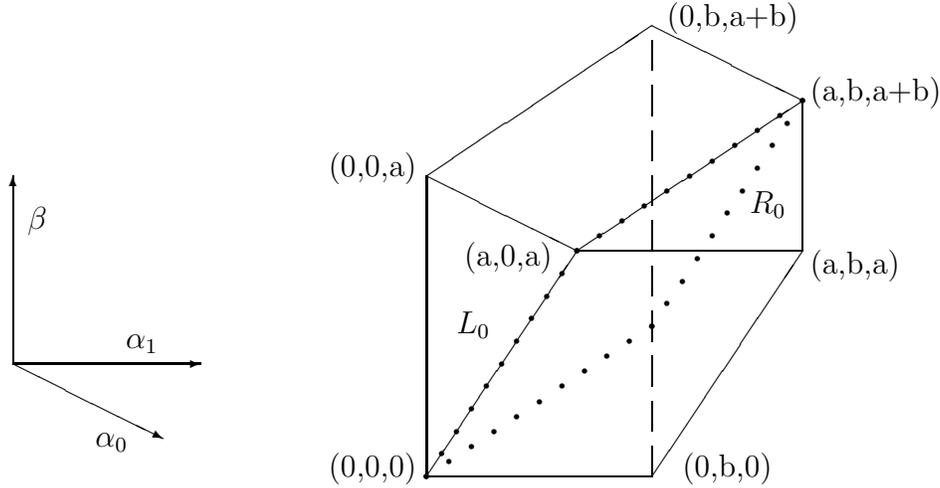
\begin{figure}[htb]                  
 \begin{center}
  \unitlength=1mm
  \begin{picture}(130,65)
\put(5,15){\vector(1,0){25}}
\put(5,15){\vector(0,1){25}}
\put(5,15){\vector(2,-1){20}}
\put(16,4){$\alpha_0$}
\put(20,17){$\alpha_1$}
\put(7,33){$\beta$}
\Xcomment{
\put(60,0){\circle*{1.0}}
\put(60,40){\circle*{1.0}}
\put(80,30){\circle*{1.0}}
\put(90,0){\circle*{1.0}}
\put(90,60){\circle*{1.0}}
\put(110,30){\circle*{1.0}}
\put(110,50){\circle*{1.0}}
}
\put(60,0){\line(1,0){30}}
\put(60,0){\line(0,1){40}}
\put(60,0){\line(2,3){20}}
\put(60,40){\line(2,-1){20}}
\put(60,40){\line(3,2){30}}
\put(80,30){\line(1,0){30}}
\put(80,30){\line(3,2){30}}
\put(90,60){\line(2,-1){20}}
\put(90,0){\line(2,3){20}}
\put(110,30){\line(0,1){20}}
\multiput(90,0)(0,6){10}{\line(0,1){4}}
\multiput(60,0)(2,3){11}{\circle*{0.6}}
\multiput(90,20)(2,3){11}{\circle*{0.6}}
\multiput(60,0)(3,2){11}{\circle*{0.6}}
\multiput(80,30)(3,2){11}{\circle*{0.6}}
\put(47,0){(0,0,0)}
\put(94,0){(0,b,0)}
\put(47,40){(0,0,a)}
\put(65,28){(a,0,a)}
\put(92,60){(0,b,a+b)}
\put(111,27){(a,b,a)}
\put(111,50){(a,b,a+b)}
\put(64,19){$L_0$}
\put(103,35){$R_0$}
  \end{picture}
 \end{center}
 \caption{The polytope $\Pscr(a,b)$. The critical section $\Pi$ is
indicated by dots.}
  \label{fig:Pscr}
  \end{figure}

So, in the nondegenerate case $a,b>0$, $\Pscr$ has 6 facets and 7
vertices. All critical vertices of $K$ are
contained in the cutting plane $\beta=\alpha_0+\alpha_1$
(cf. Corollary~\ref{cor:neq}). It intersects $\Pscr$ by the
parallelogram $\Pi$ whose vertices are $(0,0,0),(a,0,a),(a,b,a+b)$ and
a point lying on the edge of $\Pscr$ connecting $(0,b,0)$ and
$(0,b,a+b)$. This ``critical section'' $\Pi$ subdivides $\Pscr$
into two triangular prisms being, respectively, the
convex hulls of the sails in $\Lscr$ and of the sails in $\Rscr$. The
polytope $\Pscr$ is illustrated in Fig.~\ref{fig:Pscr}.

A similar description can be obtained for the convex hull $\Pscr'$ of
the RC-graph $K$ when $K$ is embedded by use of $\rho'$ in the space
$\Rset^3$ with coordinates $(\alpha,\beta_0,\beta_1)$.
Comparing~\refeq{rest} and~\refeq{rest1}, one can determine the
canonical bijection $\omega:\Pscr\to\Pscr'$ (preserving the
vertices of $K$). This $\omega$ is piecewise-linear and maps a point
$(\alpha_0,\alpha_1,\beta)$ to $(\alpha,\beta_0,\beta_1)$ such that
\begin{equation}  \label{eq:bijec}
\alpha=\alpha_0+\alpha_1,\qquad \beta_0=\min\{\alpha_1,\beta-
\alpha_0\},\qquad \beta_1=\max\{\alpha_0,\beta-\alpha_1\}.
  \end{equation}

\noindent
{\bf 4.}
Next we consider the shifted polytope $\tilde\Pscr:=
\{(0,a,0)\}+\Pscr$. For a point $(\alpha_0,\alpha_1,\beta)$ in
$\Pscr$, let $(z,x,y)$ be the corresponding point in $\tilde\Pscr$,
i.e., $z=\alpha_0$, $x=\alpha_1+a$, $y=\beta$. Following
Proposition~\ref{pr:Pscr}, $\tilde\Pscr$ is described by the linear
inequalities
  \begin{equation}  \label{eq:pattern}
0\le z\le a, \qquad a\le x\le a+b, \qquad z\le y\le x.
  \end{equation}

A triple $(x,y,z)$ (where we change the order of entries)
satisfying~\refeq{pattern} is nothing else than a real
{\em Gelfand-Tsetlin array} (or, briefly, {\em GT-array})
with border $\lambda=(\lambda_1,\lambda_2,\lambda_3)$, where
$\lambda_1:=a+b$, $\lambda_2:=a$, $\lambda_3:=0$. (The integer
triples among these correspond to the {\em semi-standard Young
tableaux} with shape $(\lambda_1,\lambda_2,\lambda_3)$; for a
definition and survey see, e.g.,~\cite{sta-99}. More about
combinatorial and polyhedral aspects of GT-arrays can be found
in~\cite{DKK-05}.)
In our case we deal with the simplest sort of such arrays, namely,
with triangular arrays of size 2; they are usually associated with
the diagram

 \begin{center}
  \unitlength=1mm
  \begin{picture}(25,16)
\put(0,0){$\lambda_1$}
\put(10,0){$\lambda_2$}
\put(20,0){$\lambda_3$}
\put(5,7){$x$}
\put(15,7){$z$}
\put(10,14){$y$}
  \end{picture}
 \end{center}

Thus, \refeq{poly} and~\refeq{pattern} explicitly indicate a
one-to-one correspondence between the set $V$ of vertices of the
RC-graph $K(a,b)$ and the set
$\Mscr$ of integer GT-arrays $(x,y,z)$ with border $(a+b,a,0)$.
According to Kashiwara~\cite{kas-95}, there is the structure of a
crystal graph on the set of GT-arrays.
(See also~\cite{ste-02,DK-05}.) For the 2-colored crystal graph
on GT-arrays with border $(a+b,a,0)$, a 1-edge (2-edge) describes
a feasible transformation in the lower row $(x,z)$ (resp. in the
upper row $(y)$) of an array in $\Mscr$. More precisely,
for $M=(x,y,z)\in\Mscr$:
  \begin{numitem}
(a) if $y<x$, then $M$ is connected by 2-edge $(M,M')$ with the array
$M':=(x,y+1,z)$; (b) if $z<a,y$ and $y-z>x-a$, then there is 1-edge
from $M$ to $(x,y,z+1)$; (c) if (b) is not applicable and if
$x<a+b$, then there is 1-edge from $M$ to $(x+1,y,z)$.
  \label{eq:p_edges}
  \end{numitem}

One can check that the edges on $\Mscr$ defined in this way
correspond to the edges of $K(a,b)$, and therefore, the 2-colored
crystal graph on GT-arrays with border $(a+b,a,0)$ is isomorphic to
$K(a,b)$. (Thus, we have a proof, alternative to~\cite{ste-03}, that
the set of (locally finite) $A_2$-regular crystals is
isomorphic to the set of RC-graphs.)

\medskip
Finally, recall that the map $\rho'$ defined in part 3 gives another
embedding of $K$ to $\Rset^3$ (with coordinates $(\alpha,\beta_0,
\beta_1)$). A point $(\alpha,\beta_0,\beta_1)$ in the polytope
$\Pscr'$ corresponds to the point $(y',z',x')$ in the shifted polytope
$\tilde\Pscr':=\{(0,0,b)\}+\Pscr'$, where $y'=\alpha$, $z'=\beta_0$,
$x'=\beta_1+b$. The corresponding analog of Proposition~\ref{pr:Pscr}
for $\Pscr'$ implies that $\tilde\Pscr'$ is described as
\begin{equation}  \label{eq:pattern1}
0\le z'\le b,\qquad b\le x'\le a+b, \qquad z'\le y'\le x',
  \end{equation}
\noindent giving the set of (real) GT-arrays $(x',y',z')$ with
border $(a+b,b,0)$. The bijection $\omega:\Pscr\to\Pscr'$ determines
a bijection $\tilde\omega$ of $\tilde\Pscr$ to
$\tilde\Pscr'$. Using~\refeq{bijec}, one can obtain an explicit
expression for $(x',y',z')= \tilde\omega(x,y,z)$:
\begin{equation}  \label{eq:bijec1}
x'=b+\max\{z,y-x+a\},\qquad y'=x+z-a,\qquad z'=\min\{x-a,y-z\}.
  \end{equation}
  %

\section{The universal RC-graph} \label{sec:univ}

By the {\em universal RC-graph} we mean the disjoint union
$\bf UC$ of RC-graphs $K(a,b)$ for all $a,b\in\Zset_+$. The
characterization of RC-graphs $K(a,b)$ given in
Section~\ref{sec:theo} and additional results from
Section~\ref{sec:patt} enable us to construct a reasonable
embedding for $\bf UC$.

In this construction, each vertex of $\bf UC$ is encoded by a tuple
$\phi=(X,x,c,D,U)\in\Zset^5$ satisfying
  \begin{equation} \label{eq:univ}
0\le x,c\le X, \quad 0\le D,U.
  \end{equation}
\noindent Moreover, there is a one-to-one correspondence between
the vertices and tuples. Under this correspondence, the
vertex set $\bf V$ of $\bf UC$ turns into a semi-group (``cone'')
$\Cscr$ in the Abelian group $\Zset^5$.

To explain the correspondence, consider a vertex $v$ in an
RC-graph $K(a,b)$, the 1-line $P_1(v)$ passing through $v$, and the
critical vertex $r$ in this line. We assign
  $$
X:=|P_1(v)|,\quad x:=\ellin_1(v),\quad c:=\ellin_1(r), \quad
D:=\ellin_2(r),\quad U:=\ellout_2(r).
  $$
Clearly~\refeq{univ} holds for these values.

Conversely, consider $\phi=(X,x,c,D,U)\in\Zset^5$
satisfying~\refeq{univ}. We associate with $\phi$ the RC-graph
$K(a,b)$, where $a:=c+U$ and $b:=X-c+D$. Then the numbers $X,c$
determine a (unique) critical vertex $r$ in $K(a,b)$, namely, $r$ is
the common vertex of the maximal left sail $L_{X-c}$ and the maximal
right sail $R_c$ (using the numeration of maximal sails as in
Section~\ref{sec:theo}). The required vertex $v$ belongs to the 1-line
passing through $c$; it is defined as having the local coordinates
$(x,c)$ in the left sail $L_{X-c}$ if $x\le c$, and coordinates
$(x-c+D,D)$ in the right sail $R_c$ if $x\ge c$.

One can explicitly express how the partial operators $F_1$ and
$F_2$ (corresponing to the 1-edges and 2-edges of $\bf UC$) act on
elements of the cone $\Cscr$. Indeed, given $\phi=(X,x,c,D,U)\in
\Cscr$, the action of $F_1$ is quite simple: $F_1$ is applicable to
$\phi$ when $x<X$, in which case it brings $\phi$ to $(X,x+1,c,D,U)$.
The operator $F_2$ is piecewise linear: (i) for $x\le c$, it is
applicable when $U>0$, in which case $F_2(\phi)=(X+1,x,c+1,D,U-1)$;
and (ii) for $x>c$, $F_2$ is always applicable and
$F_2(\phi)=(X-1,x-1,c,D+1,U)\in \Cscr$.

Identifying $\bf V$ with $\Cscr$, we observe that the trivial RC-graph
$K(0,0)$ is just the origin $\bf 0$ of the cone $\Cscr$. The RC-graph
$K(1,0)$ consists of three points $P=(0,0,0,0,1)$, $Q=(1,0,1,0,0)$,
$R=(1,1,1,0,0)$ connected by the 1-edge $(Q,R)$ and the 2-edge
$(P,Q)$. The RC-graph $K(0,1)$ consists of three points
$S=(1,0,0,0,0)$, $T=(1,1,0,0,0)$, $W=(0,0,0,1,0)$ connected by the
1-edge $(S,T)$ and the 2-edge $(T,W)$. One can check that any
nontrivial RC-graph $K(a,b)$ is obtained by taking the
Minkowsky sum of $a$ copies of $K(1,0)$ and $b$ copies of $K(0,1)$.
The cone $\Cscr$ has six ``extreme rays'',
namely, those generated by $P,Q,R,S,T,W$. The generators $P$
and $W$ are ``free'', while $Q,R,S,T$ obey the relation $Q+T=R+S$.


\end{document}